\PassOptionsToPackage{unicode}{hyperref}
\PassOptionsToPackage{hyphens}{url}
\documentclass[
  twocolumn]{article}
\usepackage{amsmath,amssymb}
\usepackage{iftex}
\ifPDFTeX
  \usepackage[T1]{fontenc}
  \usepackage[utf8]{inputenc}
  \usepackage{textcomp} 
\else 
  \usepackage{unicode-math} 
  \defaultfontfeatures{Scale=MatchLowercase}
  \defaultfontfeatures[\rmfamily]{Ligatures=TeX,Scale=1}
\fi
\usepackage{lmodern}
\ifPDFTeX\else
\fi
\IfFileExists{upquote.sty}{\usepackage{upquote}}{}
\IfFileExists{microtype.sty}{
  \usepackage[]{microtype}
  \UseMicrotypeSet[protrusion]{basicmath} 
}{}
\makeatletter
\@ifundefined{KOMAClassName}{
  \IfFileExists{parskip.sty}{%
    \usepackage{parskip}
  }{
    \setlength{\parindent}{0pt}
    \setlength{\parskip}{6pt plus 2pt minus 1pt}}
}{
  \KOMAoptions{parskip=half}}
\makeatother
\usepackage{xcolor}
\usepackage[margin=1in]{geometry}
\usepackage{graphicx}
\makeatletter
\def\maxwidth{\ifdim\Gin@nat@width>\linewidth\linewidth\else\Gin@nat@width\fi}
\def\maxheight{\ifdim\Gin@nat@height>\textheight\textheight\else\Gin@nat@height\fi}
\makeatother
\setkeys{Gin}{width=\maxwidth,height=\maxheight,keepaspectratio}
\makeatletter
\def\fps@figure{htbp}
\makeatother
\providecommand{\tightlist}{%
  \setlength{\itemsep}{0pt}\setlength{\parskip}{0pt}}
\usepackage{graphics}
\usepackage{amsmath, amsfonts, amsthm, amssymb, amscd,a4wide}
\usepackage{hyperref}
\usepackage{url}
\usepackage{float}
\usepackage{epstopdf}
\usepackage{subcaption}
\let\oldmarginpar\marginpar
\renewcommand\marginpar[1]{\-\oldmarginpar[\raggedleft\footnotesize #1]%
{\raggedright\footnotesize #1}}

\usepackage{xcolor} 
\usepackage{tikz} 
\usetikzlibrary{arrows} 
\usepackage{adjustbox}

\usepackage{slashed}
\usepackage{amsmath}
\usepackage{xcolor}
\usepackage{hyperref}
\hypersetup{linktoc = all}                
\hypersetup{hidelinks}
\hypersetup{bookmarksnumbered}
\numberwithin{equation}{section}

\newcommand \be   {\begin{equation}}
\newcommand \bel {\begin{equation}\label}
\newcommand \ee   {\end{equation}}

\newcommand \RR    {\mathbb{R}}

\newcommand \la         \langle
\newcommand \ra     \rangle

\usepackage{mathrsfs}

\usepackage{authblk}

\tikzstyle{startstop} = [rectangle, rounded corners, minimum width=3cm, minimum height=1cm,text centered, draw=black ]
\tikzstyle{io} = [rectangle, rounded corners, minimum width=3cm, minimum height=1cm,text centered, draw=black] 
\tikzstyle{process} = [rectangle, minimum width=3cm, minimum height=1cm, text centered, draw=black, fill=orange!30]
\tikzstyle{decision} = [diamond, minimum width=3cm, minimum height=1cm, text centered, text - white, draw=black, fill=black!30] 
\tikzstyle{arrow} = [thick,->,>=stealth]
\tikzstyle{bbox} = [rectangle, minimum width=3cm, minimum height=1cm, text centered, text = white, draw=black, fill=black]
\tikzstyle{pp} = [rectangle, minimum width=3cm, minimum height=0.5cm, text centered, draw=black] 
\tikzstyle{pp1} = [rectangle, draw=black!50, thick, minimum width=0.5cm, minimum height = 0.5cm]
\tikzstyle{crl} = [circle, draw=black!50, thick, minimum size = 1.5cm]
\tikzstyle{crl1} = [circle, draw=black!50, thick, minimum size = 0.7cm]
\tikzstyle{line} = [draw, -latex']
\newsavebox{\tempbox}

\tikzstyle{block} = [draw, rectangle, 
    minimum height=3em, minimum width=6em]
\tikzstyle{sum} = [draw, fill=blue!20, circle, node distance=1cm]
\tikzstyle{input} = [coordinate]
\tikzstyle{output} = [coordinate]
\tikzstyle{pinstyle} = [pin edge={to-,thin,black}]
\usepackage{booktabs}
\usepackage{longtable}
\usepackage{array}
\usepackage{multirow}
\usepackage{wrapfig}
\usepackage{float}
\usepackage{colortbl}
\usepackage{pdflscape}
\usepackage{tabu}
\usepackage{threeparttable}
\usepackage{threeparttablex}
\usepackage[normalem]{ulem}
\usepackage{makecell}
\usepackage{xcolor}
\ifLuaTeX
  \usepackage{selnolig}  
\fi
\IfFileExists{bookmark.sty}{\usepackage{bookmark}}{\usepackage{hyperref}}
\IfFileExists{xurl.sty}{\usepackage{xurl}}{} 
\hypersetup{
  pdftitle={Extrapolation and generative algorithms for three applications in finance},
  pdfauthor={Philippe G. LeFloch, Jean-Marc Mercier\^{}\textbackslash dag, and Shohruh Miryusupov},
  hidelinks,
  pdfcreator={LaTeX via pandoc}}

\title{Extrapolation and generative algorithms for
\\
three applications in finance}
\author{Philippe G.
LeFloch\footnote{Laboratoire Jacques-Louis Lions, Sorbonne University and Centre National de la Recherche Scientifique, 4 Place Jussieu, 75258 Paris, France. Email: contact@philippelefloch.org} \hskip.01cm,
Jean-Marc Mercier\(^\dag\), and Shohruh
Miryusupov\footnote{MPG-Partners, 136 Boulevard Haussmann, 75008 Paris, France. Email: jean-marc.mercier@mpg-partners.com, shohruh.miryusupov@mpg-partners.com.}}
\date{April 2024} 

\begin{document}
\maketitle
\begin{abstract}
For three applications of central interest in finance, we demonstrate
the relevance of numerical algorithms based on reproducing kernel
Hilbert space (RKHS) techniques. Three use cases are investigated.
First, we show that extrapolating from few pricer examples leads to
sufficiently accurate and computationally efficient results so that our
algorithm can serve as a pricing framework. The second use case concerns
reverse stress testing, which is formulated as an inversion function
problem and is treated here via an optimal transport technique in
combination with the notions of kernel-based \emph{encoders},
\emph{decoders}, and \emph{generators}. Third, we show that standard
techniques for time series analysis can be enhanced by using the
proposed \emph{generative algorithms}. Namely, we use our algorithm in
order to extend the validity of any given quantitative model. Our
approach allows for conditional analysis as well as for escaping the
`Gaussian world'. This latter property is illustrated here with a
portfolio investment strategy.
\end{abstract}

\hypertarget{introduction}{%
\section{Introduction}\label{introduction}}

\emph{Motivation}. We build on, and further expand, our earlier research
study \cite{LeMeI,LeMeII}, which led us to a novel class of kernel-based
algorithms and allowed us to deal with applications of central interest
in finance. Kernel-based methods are extremely efficient for financial
analytics thanks to several fundamental advantages: they provide
critical interpretability for audit and regulatory compliance, offer
robustness in sparse data scenarios, and maintain computational
efficiency which is critical for real-time analysis. Such methods have
become increasingly valuable in finance, and become recognized for their
ability to perform complex nonlinear transformations ---in turn
enhancing the predictive power of any given model without compromising
algorithmic efficiency.

\emph{State of the art and main contribution.} This Note introduces to
several applications which employ (either now standard or recently
developed) kernel-based algorithms. Among the traditional techniques, we
focus first on the so-called predictive machines, formulated as methods
for interpolation and extrapolation of data, showcasing their relevance
in different financial contexts. We present a novel generative method,
which combines ideas from optimal transport theory and kernel-based
methods. Generative methods became popular under the name of Generative
Adversarial Networks (GANs), which have emerged as a powerful class of
models and allow one to generate very realistic images. The latter have
not only captured the imagination of the public through, for instance,
Midjourney\footnote{\url{https://www.midjourney.com/home}}, 
but have also spurred extensive research
and development, eventually leading to the creation of over 500 GAN
variants\footnote{\url{https://github.com/hindupuravinash/the-gan-zoo}}.
%
 Among the diverse families of GANs, conditional
GANs stand out for their unsupervised learning capability and produce
images based on conditional inputs, while TimeGan \cite{yoon2019} allows
the generation of realistic time series data. In parallel with the
evolution of GANs, another significant strand of research in generative
modelling takes it roots in traditional kernel methods, such as Sinkhorn
autoencoders \cite{patrini2018} and Schr\"oodinger bridge generative
models \cite{Pham:2023} for image and time series data. In agreement
with optimal transport-based generative techniques, our methodology for
constructing suitable algorithms was inspired by the method of
Nadaraya-Watson kernel regression \cite{Nadaraya:1964}, which guided us
to design our sampling algorithm for conditioned density estimations. At
the core of our method is a sampling algorithm that employs optimal
transport and maps a white-noise latent space directly to the target
distribution space. This process is uniquely implemented by using
permutation indices of samples: the mappings arising in the computation
can be traced back and interpreted in the applications, which is an
essential feature in complex distribution analysis.

\emph{Applications in finance.} We illustrate the performance of our
algorithms in three practical applications within quantitative finance.

\begin{itemize}
\item \textbf{Online predictions of PnL and its sensitivities.} By adapting our prediction
algorithm, we provide a framework for online forecasting of profit and loss (PnL) statements and their sensitivities, enhancing the robustness of financial decision-making processes.

\item \textbf{Reverse-stress test (PnL function inversion).} A novel application of our permutation algorithms allows us to perform reverse-stress testing. By inverting the PnL function, we can better understand the conditions that lead to extreme financial outcomes, preparing for worst-case scenarios in risk management.

\item \textbf{Quantitative models and generative methods.}  Traditional financial models frequently rely on predefined stochastic processes, such as Brownian motions, which might not capture the full spectrum of market complexities and dynamics. Our methodology progresses in two significant directions: Initially, we ascertain that a majority of quantitative models can be conceptualized as mappings, transforming time series data into white noise. This insight paves the way for our novel contribution, where we employ optimal permutations and mappings to enhance the sophistication of quantitative models, exemplified through the refinement of the GARCH process. This innovation marks a significant advancement in our ability to model intricate market behaviors. Additionally, we employ a conditional probability estimator to devise a portfolio management strategy, anchored in precise market indicators. This strategy not only deepens our understanding of market dynamics but also facilitates the formulation of superior investment strategies.
\end{itemize}

\hypertarget{performing-stress-tests-via-an-extrapolation-algorithm}{%
\section{Performing stress tests via an extrapolation
algorithm}\label{performing-stress-tests-via-an-extrapolation-algorithm}}

\emph{Aim.} We present our extrapolation algorithm, particularly
designed for asset pricing applications when dealing with potentially
large multi-asset portfolios. Our motivation for proposing this test
stems from the fact that pricing engines in the finance industry are
often computationally intensive and therefore struggle with real-time
deployment. Our aim here is to demonstrate that learning a pricing
function offline and then using extrapolation provide an approach that
is sufficiently accurate for dynamical use as a real-time risk/pricing
framework for stress testing or for Profit and Loss (PnL) analysis.

\begin{figure}

{\centering \includegraphics[width=0.9\linewidth]{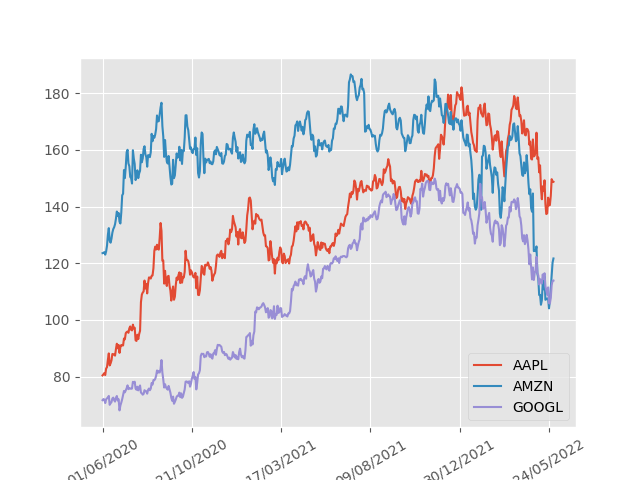} 

}

\caption{\label{fig:AAG}Charts for Apple, Amazon, and Google}\label{fig:unnamed-chunk-3}
\end{figure}

\emph{Notation: market data and time series.} We consider time series
denoted by \(t \mapsto X(t)\in \mathbb{R}^D\) which is observed on a
time grid \(t^{-T_x}<\ldots<t^{0}\). In this context, \(t^{0}\) denotes
the pricing date and our notation is as follows:
\begin{equation}\label{TS}
X=\Big(x^{n,k}_d\Big)_{d=1\ldots D}^{k = 1,\ldots,T_x} \in \mathbb{R}^{D,T_x}.
\end{equation} In our example, the observed data comprise \(253\)
closing values denoted \(x^{-252},\ldots,x^{0}\), for the S\&P500 market
index during the period of time from \(t^{-252}=\) June 1, 2021 to
\(t^{0}=\) June 1, 2022. (These data were retrieved from Yahoo Finance.)
In the notation \eqref{TS}, this dataset is represented by a matrix with
dimensions \(D=3, T_x = 253\) and the corresponding charts are displayed
in Figure \ref{fig:AAG}.

The integer \(T_x\) is the number of historical observations of the time
series \(X\), \(D\) represents the number of components of the observed
process, and, for the sake of simplicity in the presentation, our
example is based on three assets.

\begin{figure}

{\centering \includegraphics[width=1\linewidth]{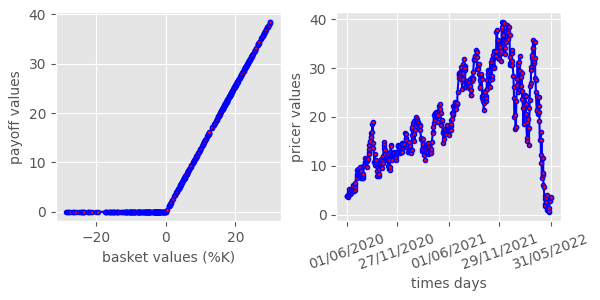} 

}

\caption{\label{fig:pricer}Payoff (left), and pricer (right) values}\label{fig:unnamed-chunk-4}
\end{figure}

\emph{Notation: portfolio.} We focus on a function
\(P(t,x) \in \mathbb{R}^{D_P}\), where \(x \in \mathbb{R}^D\) and which
represents an external pricing engine, evaluating a portfolio of \(D_P\)
instruments, based on assets valued at \(x\) at time \(t\). This setup
allows us to calculate the portfolio's value at any given time,
especially at its maturity \(T\), where the payoff is specifically
defined as \(P(T,x)\). To simplify our presentation, we assume a single
instrument, \(D_P=1\), focusing on a basket option with underlying
assets. The payoff for this basket option is computed as
\(P(T,x) = \max(\langle \omega, x \rangle - K, 0)\), where
\(\langle \omega, x \rangle\) represents the weighted sum of the
basket's asset values with \(\omega\) being the equal weights. Here,
\(K\) represents the strike price of the option and \(T\) denotes its
maturity. Figure \ref{fig:pricer} displays, on the left-hand side, the
payoff as a function of the basket values and, on the right-hand side,
the pricing engine values as a function of time. For demonstration
purposes, we use a simplified Black-Scholes model, denoted as
\(P(t,x) := \text{BS}(S, T-t, \sigma, K)\) with \(S=K=<\omega,X^0>\),
implying the initial basket value equals the strike price. Of course,
our framework can accommodate many different types of pricing functions.

\emph{A kernel-based extrapolation algorithm.} Using historical data, we
generate synthetic data for a future date \(t^1 = t^0 + H\) days through
time series forecasting, as described earlier, with \(H=10\) days in
this example, in order to simulate stress tests scenarios. These data
are denoted by
\(X := \big(x^n_d\big)_{d=1\ldots D}^{n=1\ldots N_x} \in \mathbb{R}^{N_x,D}\).
Similarly, we produce another set of scenarios
\(Z \in \mathbb{R}^{N_z,D}\) using the same methodology. Consequently,
we simulate learning the pricing engine values \(P(X) \equiv P(t^1,X)\)
from the market data \(X\), referred to as the \emph{training set}, and
then extrapolate this function onto the set \(Z\), referred to as the
\emph{test set}, to compare against the actual pricing engine values
\(P(Z)\).

\begin{figure}

{\centering \includegraphics[width=1\linewidth]{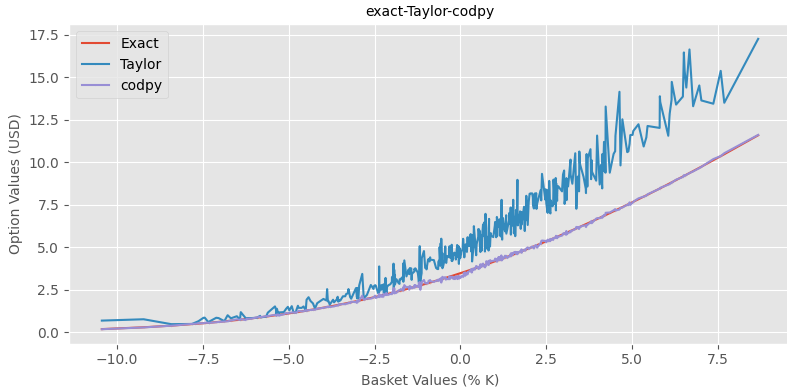} 

}

\caption{\label{fig:output}A benchmark of PnL extrapolation methods}\label{fig:unnamed-chunk-5}
\end{figure}

To describe the extrapolation procedure, we introduce basic notions
relative to kernel operators, and refer to
\cite{BerlinetThomasAgnan:2004} for a comprehensive overview of
reproducing kernel Hilbert space RKHS theory. Let
\(\mathcal{X} \subset \RR^D\) a convex set, we call a function
\(k: \mathcal{X} \times \mathcal{X} \mapsto \mathbb{R}\) a kernel if it
is a symmetric and positive definite (see
\cite{BerlinetThomasAgnan:2004} for a definition). If
\((X,Y):=\Big(x^{n}_d\Big)_{d=1\ldots D}^{n=1\ldots N_x},\Big(y^{n}_d\Big)_{d=1\ldots D}^{n=1\ldots N_y}\)
are two set of distinct points, we define
\(k(X,Y):=k(x^i,y^j)_{i=1 \ldots N_x}^{j=1 \ldots N_y}\) the Gram
matrix. Let \(P\) be any function taking values on \(\mathcal{X}\), and
denote \(P(X):=P(x^1),\ldots,P(x^{N_X})\) its discrete values. Then, for
any \(z \in \RR^D\), we define the \emph{projection} as
\begin{equation} 
\label{projection}
P_k(X,Y)(z) := k(Y,z) k(X,Y)^{-1}P(X).
\end{equation} In this equation, the inverse of the Gram matrix
\(k(X,Y)^{-1}\) is computed by applying a standard least-square method.
Extrapolation is defined as \(z\mapsto P_k(X,X)(z)\), denoted
\(P_k^X(z)\) for short, and is a \emph{reproducible} operation, in the
sense that it satisfies \(P_k(X) = P(X)\); namely, it is exact on the
training set. In a similar way, we can define other types of derivative
operators such as, for instance, the gradient operator \begin{equation} 
\label{nabla}
\nabla P_k(X,Y)(z) := (\nabla k(Y,z)) k(X,Y)^{-1}P(X),
\end{equation} and we denote \(\nabla P_k^X(z)\) for short in the
extrapolation mode. To benchmark the proposed approach, we compute the
prices on the test set \(Z\) using three methods and display our results
in Figure \ref{fig:output}.

\begin{itemize}
\item Analytical prices : it is computed as $P(t^1,Z)$, and is the reference values for our tests.

\item Predicted prices or PnL: the price function $P(t^1,Z)$ is approximated using the formula \eqref{projection}, as $\mathcal{P}_{k}(Z)$.

\item $\Delta$-$\Gamma$ approximation: the price function $P(t^1,Z)$ is computed using a second order Taylor formula approximation around $P(t^0,x^0)$, hence involving the price, four derivatives (three for each asset, one for time), and sixteen second order derivatives.

\footnote{We benchmarked against a Taylor approximation, as this method is currently used by some banks to estimate their PnL on a real time basis.}
\end{itemize}

\begin{figure}

{\centering \includegraphics[width=1\linewidth]{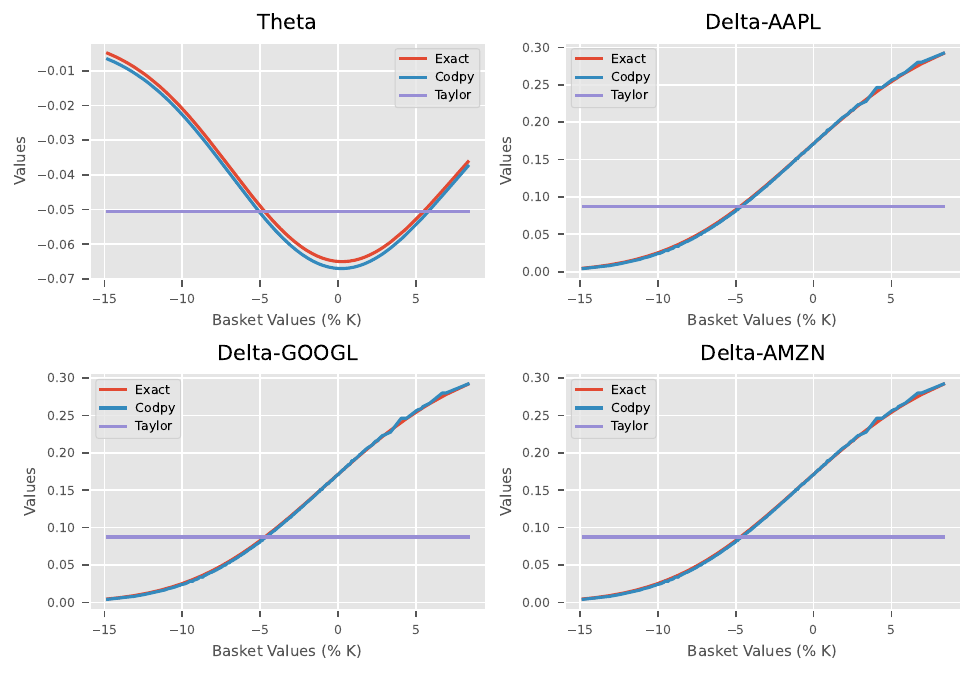} 

}

\caption{\label{fig:greeks}A benchmark of PnL greeks from an extrapolation methods}\label{fig:unnamed-chunk-6}
\end{figure}

In the same vein, we can also compute greeks using \eqref{nabla},
resulting in four plot, three deltas for each asset, one theta for time,
in Figure \ref{fig:greeks}.

\emph{In conclusion.} The above experiments show that extrapolation
methods achieve basis point accurate
method\footnote{See \cite{LeMeMi} for further testing.}, even with a
limited number of pricing examples, demonstrating their viability for
real-time pricing applications. We point out that our approach ofers
multiple ways for further improvement and optimization. For instance,
the selection of stress test scenarios, denoted by \(X\), can be refined
through clustering-based strategies, for instance by relying on
\emph{sharp discrepancy} sequences; see \cite{LeMeI}.

\hypertarget{reverse-stress-tests-as-encoders}{%
\section{Reverse stress tests as
encoders}\label{reverse-stress-tests-as-encoders}}

\emph{Aim.} Reverse stress testing stands in contrast to traditional
stress testing by starting with a specific outcome, such as portfolio
values or PnL losses denoted as \(p \in \mathbb{R}^{D_p}\), and
backtracking to uncover the market scenarios that could lead to such
outcomes. This approach is invaluable for identifying vulnerabilities
within a portfolio and enhancing risk management strategies. We
illustrate this concept by considering portfolio or PnL values
\(p \in \RR^{D_p}\), aiming to reverse-engineer the market data
scenarios \(x=P^{-1}(t^1,p)\) using a pricer function \(P\), known from
discrete values \(P:=P^1,\ldots,P^{N_X} \in \RR^{N_x,D_p}\), where
\(P^n=P(t^1,X^n)\) is evaluated on market data
\(X:= X^1,\ldots,X^{N_x}\). The challenge arises when this mapping,
\(x \mapsto P(x)\), from \(\mathbb{R}^D\) to \(\mathbb{R}^{D_P}\), lacks
an obvious inverse due to its non-invertibility. We address this
challenge by introducing the concept of encoders and decoders, borrowed
from machine learning, but using kernel-RKHS methods.

\emph{Encoders, decoders and generators.} Encoders in this particular
numerical experiment are conceptualized as smooth, invertible maps,
\(x \mapsto P(x)\) from \(\RR^D\) to \(\RR^{D_P}\) that bridge the gap
between the market data \(X\) and the portfolio values \(P\), and we
denote this maps as \(x\mapsto \mathcal{L}(X,P)(x)\), while its inverse
is formally denoted \(p\mapsto \mathcal{L}(P,X)(p)\). The extrapolation
operator \eqref{projection}, denoted by the equation below, initially
suggests a direct inversion approach: \begin{equation} \label{NRST}
 k(P,p) k(P,P)^{-1}X.
\end{equation} Yet, this direct approach may falter due to the inherent
non-invertibility of \(P\) as a mapping from \(\RR^{D_p}\) to
\(\RR^{D}\). To enhance stability, we propose a refined method, relying
on a permutation \(\sigma:[1,\ldots,N_x] \mapsto [1,\ldots,N_x]\) of the
original data set, written as \begin{equation} \label{NRST=two}
\mathcal{L}(P,X)(p) := k(P,p) k(P,P)^{-1}X^\sigma
\end{equation}

In the particular case where spaces have matching dimensions, that is
\(D_x=D\) in our notation (cf.~also \cite{Brezis:2018}), considering a
distance function \(d(x,y)\), optimal transport theory proposes to
determine this permutation as \begin{equation} \label{OT}
\overline{\sigma} = \arg \inf_{\sigma \in \Sigma} \sum_{n=1}^{N_x} d(X^{\sigma(n)},P^n)
\end{equation} where \(\Sigma\), the set of all permutations. For
kernels methods, a natural distance is given by
\(d_k(x,y)=k(x,x)+k(y,y)-2k(x,y)\), called the kernel discrepancy or
maximum mean discrepancy, see \cite{GR:2006}.

However, the existing literature seems less profuse when the function
maps unrelated spaces, that is \(D_x \neq D\) in our notation. Hence we
introduced the following Ansatz in this case, based on the gradient
formula \eqref{nabla} \begin{equation} \label{SP}
\overline{\sigma} = \arg \inf_{\sigma \in \Sigma} \| (\nabla k(P,P)) k(P,P)^{-1}X^{\sigma} \|_2^2.
\end{equation}

This approach is reminiscent of a generalized
\textit{traveling salesman problem}\footnote{see wikipedia page \url{https://en.wikipedia.org/wiki/Travelling_salesman_problem} }
and is not equivalent to the one in \ref{OT}, however it allows for the
determination of a smooth mapping between the original and target
spaces, as desired.

Encoders and decoders allows to define the notion of
\textbf{generators}, needed later on. Let \(\mathbb{X},\mathbb{Y}\) two
continuous distributions taking values in \(\RR^{D_x},\RR^{D_y}\), and
consider \(X \in \RR^{N,D_x}, Y \in \RR^{N,D_y}\) two variates of equals
length. Then \eqref{NRST} provides a method to \emph{generate} a sample
\(y\), from a sample \(x\), according to the formula
\(y = k(X,x) k(X,X)^{-1}Y^\sigma\). In particular, if one consider
\(\mathbb{X}\) as a known distribution, as for instance a uniform one,
this provides a modeling of any distribution \(\mathbb{Y}\) by a
continuous one, that is statistically similar to the observed data
\(Y\).

\emph{In conclusion.} Using \eqref{SP}, we produced the following
picture, corresponding to a reverse stress test, that we comment
thereafter.

\begin{figure}

{\centering \includegraphics[width=1\linewidth]{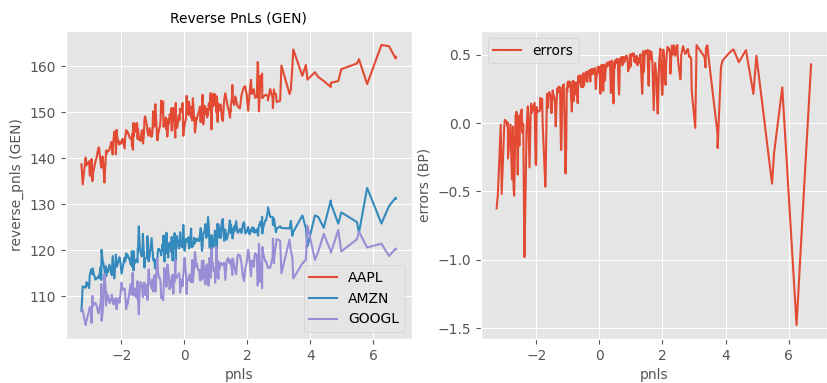} 

}

\caption{\label{fig:ReversePrices} Reverse prices (left) and benchmark (right)}\label{fig:unnamed-chunk-7}
\end{figure}

Consider the original set of market data \(X=X^1,\ldots,X^{N_x}\) and
corresponding PnL \(P=P^1,\ldots,P^{N_x}\). Using a generator, as
described in the previous section, we sampled new examples
\(\overline{P}=\overline{P^1},\ldots,\overline{P^{N_x}}\). These
simulated distribution of prices is used with the reverse stress test
method to compute the corresponding scenarios at left, notation
\(\overline{X}=\overline{X^1},\ldots,\overline{X^{N_x}}\). To benchmark
this result, we computed the error distribution
\(P(\overline{X^1})-\overline{P^1}),\ldots,P(\overline{X^{N_x}}-\overline{P^{N_X}})\).
This distribution is plot as the red line at right, expressed in basis
points.

\hypertarget{modeling-time-series-via-a-generative-algorithm}{%
\section{Modeling time series via a generative
algorithm}\label{modeling-time-series-via-a-generative-algorithm}}

\emph{Mapping time series.} Our approach consists in a general model
framework where a time series \(X\) is transformed into a latent
variable \(\varepsilon\), interpreted as white noise, through a
continuous and invertible map \(F\). With this approach, the
transformation \(F(X) = \varepsilon\) isolate the inherent noise within
the data, facilitating a deeper analysis and manipulation of the
underlying stochastic processes. Let us give a first simple example to
fix ideas. Consider any time series \(X\) of the form \eqref{TS}, then
the following simple random walk provides such a modeling: \[
\varepsilon^k = X^{k+1}-X^k, F^{-1}(\varepsilon)^k=X^0+\sum_{n\le k} \varepsilon^n.
\] A time series model can then be interpreted as an invertible mapping
from \(X \in \RR^{D_x,T_x}\) into
\(\varepsilon \in \RR^{D_{\varepsilon},T_{\varepsilon}}\). As
demonstrated wiht the GARCH model in the following, most quantitative
models can be interpreted through such mappings \(F\), and we provide
numerous other examples in \cite{LeMeMi}. The purpose of such an
approach is to provide a framework to deal with time series enjoying the
following properties.

\begin{itemize}
\tightlist
\item
  The time series \(X\) is \emph{reproducible}, meaning that
  \(X=F^{-1}(\varepsilon)\). This property ensures that \(X\), the
  historical dataset, is in the range of the model.
\item
  The variable \(\varepsilon \in \RR^{D_{\varepsilon},T_{\varepsilon}}\)
  can be reproduced from any other random variable
  \(\eta \in \RR^{D_{\eta},T_{\eta}}\), for instance a uniform sampling,
  using a \emph{generator}, as described earlier. Note that the
  distribution that generates the variate \(\varepsilon\) can be rather
  arbitrary, and we are not restricted to assuming Gaussian
  distributions for time series, as is usually the case in classical
  modeling of time series.
\item
  It provides a simple Monte-Carlo framework, as sampling new
  \(\overline{\varepsilon}\) generates a new trajectory defined as
  \(\overline{X} = F^{-1}(\overline{\varepsilon})\).
\end{itemize}

From an economical point of view, these models simulates new
trajectories that are calibrated to the observed dynamics of the
financial time series, as they rely on reproducing historically observed
white noise. This way of modeling ensures that short-term dynamics of
time series are easily and better captured. It provides also a powerful
tool for forecasting, risk assessment, and strategy benchmarking, as we
can apply this framework in the following areas.

\begin{itemize}
\tightlist
\item
  Benchmarking strategies, where we resample the original signal \(X\)
  on the same time-lattice to draw several simulated trajectories and
  compare them to the original one using various performance indicators.
\item
  Monte-Carlo forecast simulations for future time points, allowing for
  the exploration of potential future scenarios based on historical
  data.
\item
  Forward Calibration, where we frame it as a minimization problem with
  constraints, optimizing the generation of new samples that meet
  specific financial criteria.
\item
  PDE pricers, enabling the computation of forward prices or
  sensitivities by solving backward Kolmogorov equations in a
  multidimensional tree structure, see \cite{LeMeI}.
\end{itemize}

\emph{The GARCH(p,q) model as an example.} We show an example of such an
model extension considering the GARCH model, in order to illustrate its
flexibility and adaptability to any quantitative model, but also its
potential to enhance the modeling and simulation of financial time
series. This approach is not exclusive to GARCH models; it is equally
applicable to the broader ARIMA family, as well as a wide range of both
discrete and continuous stochastic processes, as local or rough
volatility ones.

The generalized autoregressive conditional heteroskedasticity (GARCH)
model, particularly in its \((p, q)\) order formulation, captures
financial markets' volatility clustering---a phenomenon where
high-volatility events tend to cluster together in time. The defining
equations of a GARCH\((p,q)\) model are as follows: \[
\begin{cases}
X^k = \mu + \sigma^k Z^k, \\
(\sigma^k)^2 = \alpha_0 + \sum_{i=1}^p \alpha_i (X^{k-i})^2 + \sum_{i=1}^q \beta_i (\sigma^{k-i})^2,
\end{cases}
\]

where \(\mu\) is the mean, \(\sigma^k\) denotes the time-varying
volatility, and \(Z^k\) symbolizes a white noise process. The
coefficients \(\alpha_i\) and \(\beta_i\) govern the model's
responsiveness to changes in volatility and the inertia of past
volatility, respectively.

The variance equation of GARCH can be compactly expressed using the
backshift operator \(B\), leading to \[
(1-\beta(B))(\sigma^k)^2 = \alpha_0 + \alpha(B)(X^k)^2,
\] where \(\alpha(B) = \sum_{i=1}^p \alpha_i B^i\) and
\(\beta(B)=\sum_{i=1}^p \alpha_i B^i\) and backshift
\(B^iX^k = X^{k-i}\) allows for the compact expression of lagged effects
on the stochastic process. By defining
\(\varphi(B) = \alpha_0 + \sum_{i=1}^p \alpha_i B^i\),
\(\theta(B) = 1 - \sum_{i=1}^q \beta_i B^i\), the stochastic variance
\(\sigma^k\) can be derived as \[
\sigma^k = \sqrt{\varphi^{-1}(B) \theta(B)(X^k)^2} = \sqrt{\pi(B)(X^k)^2}.
\] with \(\pi(B)=\varphi^{-1}(B) \theta(B)\) encapsulating the variance
transformation.

Our methodology integrates the GARCH model within a broader generative
framework, identifying the \emph{`GARCH map'}, \(G: X^k \mapsto Z^k\),
which translates the observable time series \(X^k\) into a latent white
noise process \(Z^k\): \[
Z^k = G(X^k) = \sqrt{\pi^{-1}(B)(X^k)^2}(X^k - \mu).
\] This GARCH map forms the basis of our three-stage generative process
for financial time series, specifically demonstrated through the example
of stock prices in Figure \ref{fig:AAG}.

\begin{itemize}
\item
  By applying the GARCH map to the historical stock prices \(X\), we
  extract the latent noise \(\varepsilon = G(X)\). This latent noise
  encapsulates the fundamental stochastic processes driving the
  volatility modeled by the GARCH framework.
\item
  Utilizing the permutation algorithm \eqref{OT} as a generator, we
  produce new instances of the latent variable \(\tilde{\varepsilon}\),
  simulating alternative realizations of the market's stochastic
  behavior.
\item
  Employing the inverse GARCH transformation, symbolically represented
  as \(G^{-1}\), we map these new latent samples back to the domain of
  stock prices, thereby generating new trajectories
  \(\tilde{X} = G^{-1}(\tilde{\varepsilon})\), consistent with the GARCH
  dynamics.
\end{itemize}

\begin{figure}
\includegraphics{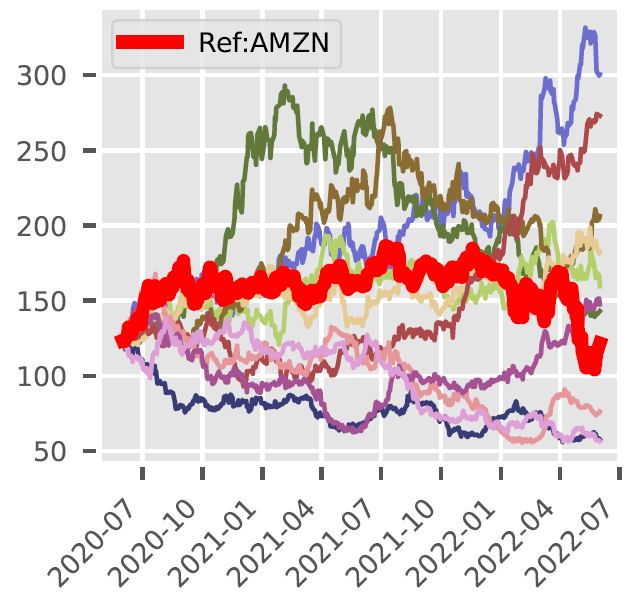}
\caption{Ten examples of generated paths with the GARCH(1,1) model} \label{fig:plot362}
\end{figure}

Figure \ref{fig:plot362} displays ten simulated trajectories of one of
the three stock, Amazon's ones, using the GARCH(1,1) model, the others
two producing quite similar patterns. This visualization exemplifies how
the generative process can give us a diverse set of plausible future
paths for the stock price, rooted in the historical volatility patterns
captured by the GARCH model.

\hypertarget{extending-quantitative-models-through-conditioning}{%
\section{Extending quantitative models through
conditioning}\label{extending-quantitative-models-through-conditioning}}

\begin{figure}

{\centering \includegraphics[width=1\linewidth]{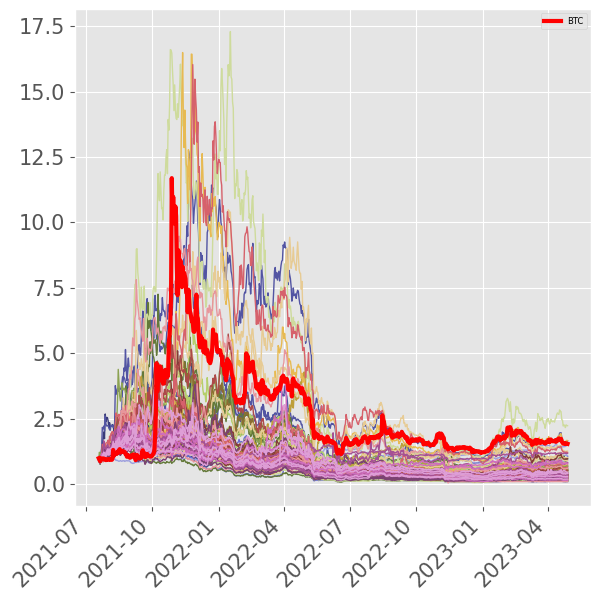} 

}

\caption{\label{fig:cryptos} A time series of 106 crypto assets}\label{fig:unnamed-chunk-8}
\end{figure}

\emph{Further applications.} We finally demonstrate that the generative
algorithms presented so far can be safely used for conditional analysis;
this is done on a numerical experiment in the context of portfolio
management. For the sake of illustration, we consider a time series
given by closing of a basket of 106 crypto currencies during the period
from 19/07/2021 to 28/04/2023, which was observed on a daily basis and
corresponds to a time series \(X\) with \(D,T_x = 106, 649\). We rely
here on our notation \eqref{TS}, and we present the plot in Figure
\ref{fig:cryptos} after suitable normalization at the initial time.

To this aim, consider a integer \(W\) defining a sliding window
\(t^{n+k},k = 0,\ldots,W\), \(n=0,\ldots T_x-W\). At each time labelled
\(n\), we define several portfolios determined by their components
\(\omega_{n,d}\), and their wealth
\(P^n := <\omega_n, X^{n+W}> = \sum_d \omega_{n,d} x_d^{n+W}\). Each
strategy corresponds to an investment at the time \(t^{n+W}\), and
portfolio performances are compared together with, and without,
conditional analysis. In the next paragraph, we describe the investment
strategies that we used, and how conditional analysis can be used to
enhance them.

Portfolio management usually build upon a simple quantitative model,
determined by \emph{returns} of the assets, that is, following the
notation introduced in the previous section, returns are determined as
the discrete distributions
\(\epsilon_n = \Big( \frac{X^{n+k+1}}{X^{n+k}}-1\Big)_{k=1,\ldots,W}\),
and the inverse map is
\(F_n^{-1}(\epsilon)= \Big( X^{n} \ \Pi_{n \le k}(1+\epsilon^{n+k}) \Big)_{k=0,\ldots}\).

\emph{Efficient portfolio.} Let us recall the Markowitz mean-variance
optimization, providing a classical investment strategy. It consists in
finding a set of portfolio weights \(\omega:=(\omega_d)_{d=1,\ldots,D}\)
as the solution of the quadratic programming problem
\begin{equation} \label{EP}
  \overline{\omega}(\epsilon) = \arg \inf_{\omega} \frac{1}{2}\omega^T Q \omega-\epsilon \omega^T \overline{\epsilon} + \beta(|\omega-\omega^0|). 
\end{equation} The terms \(\beta\) represent the transaction costs,
which are proportional to the portfolio change in weights, symbolized
here in the term \(|\omega-\omega^0|\). This comes usually with
constraints over the weights, and we used \(\sum_d \omega_d=0\) (long /
short strategy) together with \(|\omega_i| \le 1\).

As our strategies are determined via a random variable modeling the
asset returns \(\epsilon\)
\footnote{as point out earlier, this conditional method can be used with any other quantitative modeling of the assets},
we consider two quite comparable situations, where this random variable
is conditioned by two different observed random variable \(\eta\), that
is we considered the distribution \(\epsilon | \eta = \eta^{n}\) for
each time \(t^n\), \(n=W,\ldots T_x\). Now, we are going to describe how
the generative methods of the previous section can be applied to
conditioning random variables.

\emph{Conditioned random variables and latent spaces.} Let
\(\mathbb{X} \in \mathbb{R}^{D_X},\mathbb{Y} \in \mathbb{R}^{D_Y}\) two
dependent random variables and consider
\(\mathbb{Z} = \Big(\mathbb{X}, \ \mathbb{Y}\Big) \in \mathbb{R}^{D_X+D_Y}\)
the joint random variable. Let
\(Z = (X,Y) \in \mathbb{R}^{N,D_x + D_y}\) be a variate of
\(\mathbb{Z}\). Relying on \eqref{NRST}, we consider another variate
\(\epsilon\) drawn from any known distribution, say
\(\epsilon = (\epsilon^n)_{n=1\ldots N}\), which is called the
\emph{latent} and is decomposed as
\(\epsilon:=(\epsilon_x,\epsilon_y)\). Elaborating on the two encoders
map in \eqref{NRST}, namely \(\mathcal{L}(X,\epsilon_x)\) and
\(\mathcal{L}(\epsilon,Z)\), the following construction provides a
generator of the conditioned law \(\mathbb{Y} | \mathbb{X}=x\): \[
  \eta_y \mapsto \mathcal{L}(\epsilon,Y)(\eta_x,\eta_y), \quad \eta_x = \mathcal{L}(X,\epsilon_x)(x). 
\] There is a lot of freedom in choosing the distributions
\(\epsilon_x,\epsilon_y\). We can pick up, for instance, uniform
distributions, or the trivial map \(\epsilon_x = X\), which can be handy
in certain applications.

\begin{figure}
\includegraphics{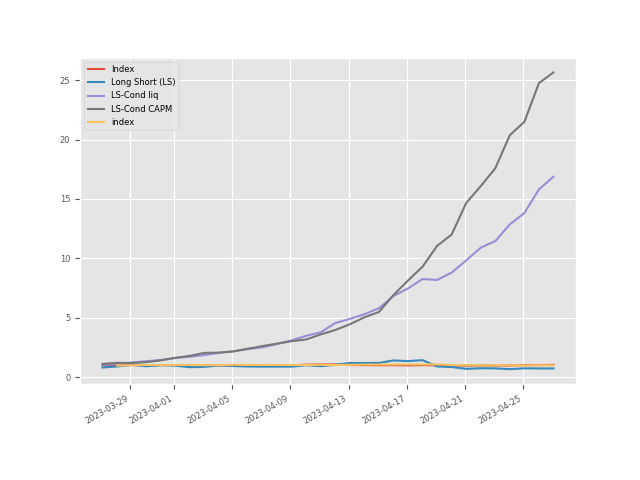}
\caption{Benchmark of different portfolio strategies} \label{fig:plot363}
\end{figure}

\emph{Benchmarks results.} Next, we provide an illustration of
allocation strategies. We compared four strategies in figure
\ref{fig:plot363}, listed as follows.

\begin{itemize}
\tightlist
\item
  The first one is given by an equiweighted portfolio, which we call
  `index' and serves as a reference.
\item
  The second one (LS) is a Long Short strategy, determined as
  approximating the mean variance problem \eqref{EP}.
\item
  The third one is also a Long Short strategy, but where the return
  distribution at time \(t^n\), \(\epsilon^{n}\), is conditioned to the
  capital asset pricing model, that is,
  \(b^n=r_f+\beta^k(\overline{\omega^k}-r_f)\), in which \(\beta^k\) is
  the regression coefficients of the weights \(\omega^k\).
\item
  The third one is very similar to the previous one, except that the
  distribution \(\epsilon^{n}\) at time \(t^{n}\) is conditioned with a
  distribution \(b^n\) contains for each assets liquidity values, time
  averages on different windows as well as differences between them.
\end{itemize}

We reported the daily performance of each portfolio
\(P_j\),\(j=1,..,4\), as
\(P_j^{n}:= \Pi_{k \le n} \frac{<\omega_j^{k},X^{k+1}>}{<\omega_j^{k},X^{k}>}\).

As we can see, strategies based on conditional returns outperformed non
conditioned ones on these examples. Of course, this experiment has to be
considered as a simple validation of our algorithms, as we conditioned
our returns values \(\epsilon^n\) with distributions that are known at
time \(t^n\). In an operational context, the available information for
conditioning comes from the past, that is the time \(t^{n-1}\). In
particular, given a conditioning distribution, this test does not answer
the quite interesting question whether there exists a time frequency of
observations below which such an approach could be profitable or not.


\begin{thebibliography}{9}

\bibitem{BerlinetThomasAgnan:2004} 
{\sc A. Berlinet and C. Thomas-Agnan,} 
{\sl Reproducing kernel Hilbert spaces in probability and statistics,} 
Springer US, Kluwer Academic Publishers, 2004.


\bibitem{Brezis:2018} 
{\sc H. Brezis,}
Remarques sur le probl\`eme de Monge–Kantorovich dans le cas discret, 
Comptes Rendus Math. 356 (2018), 207--213.


\bibitem{GR:2006}
{\sc A. Gretton, K.M. Borgwardt, M. Rasch, B. Sch\"{o}lkopf, and A.J. Smola,}
A kernel method for the two sample problems, 
Proc. 19th Int. Conf. on Neural Information Processing Systems, 2006, pp.~513--520. 

\bibitem{Pham:2023} {\sc M. Hamdouche, P. Henry-Labordere, and H. Pham,} 
Generative modeling for time series via Schr\"odinger bridge. Available as 
ArXiv:2304.05093. 

\bibitem{LeMeI}  {\sc P.G. LeFloch and J.-M. Mercier,} 
A class of mesh-free algorithms for some problems arising in finance and machine learning, 
J. Scientific Comput. 95 (2023), 75.

\bibitem{LeMeII}  {\sc P.G. LeFloch and J.-M. Mercier,} 
The transport-based mesh-free method: a short review, 
Wilmott journal, 2020, pp.~52--57. 

\bibitem{LeMeMi} {\sc P.G. LeFloch, J.-M. Mercier, and S. Miryusupov,} 
{\sl CodPy: a Python library for numerics, machine learning, and statistics,} 
Monograph, 2024. Available as ArXiv.2402.07084. 

\bibitem{liu2018} {\sc H.D. Liu, Y. Guo, N. Lei, Z. Shu, S.-T. Yau, D. Samaras, and X. Gu,} 
Latent space optimal transport for generative models, 2018. Available as ArXiv:1809.05964.

\bibitem{Nadaraya:1964}
{\sc E. A. Nadaraya,}
On estimating regression,
Theory Proba. and Appl.. 9 (1964)), 141. 

\bibitem{patrini2018} {\sc G. Patrini, R. van den Berg, P. Forré, M. Carioni, S. Bhargav, M. Welling, T. Genewein, and F. Nielsen,} 
Sinkhorn auto-encoders, 2018. 
Available as ArXiv:1810.01118. 



\bibitem{yoon2019} {\sc J. Yoon, D. Jarrett, and M. van der Schaar,} 
Time-series generative adversarial networks, 
Neural Information Processing Systems (NeurIPS), 2019.  

\bibitem{zhang2019} {\sc O. Zhang, R.-S. Lin, and Y. Gou,}
Optimal transport-based generative autoencoders, 2019. Available as ArXiv:1910.07636.

\end{thebibliography}
\end{document}